\documentclass[secthm,seceqn,amsthm,ussrhead,12pt]{amsart}
\usepackage{amsmath,latexsym}
\usepackage[english]{babel}
\usepackage[psamsfonts]{amssymb}
\usepackage{times}
\usepackage{cite}
\usepackage{pdflscape} 
\usepackage{ulem}
\usepackage[mathcal]{euscript}
\usepackage{tikz}
\usepackage{hyperref}
\usepackage{cancel}
\usetikzlibrary{arrows}

\newtheorem{Th}{Theorem}
\newtheorem{Lem}[Th]{Lemma}

\newtheorem{corollary}[Th]{Corollary}

\newenvironment{Proof}[1][Proof.]{\begin{trivlist}
\item[\hskip \labelsep {\bfseries #1}]}{\flushright
$\Box$\end{trivlist}}

\begin{document}
	\sloppy


{\Large Degenerations of Leibniz and anticommutative algebras 
\footnote{The work was supported by 
FAPESP 17/15437-6, 17/21429-6; 
AP05131123 "Cohomological and structural problems of non-associative algebras";
RFBR 18-31-00001;
the President's Program "Support of Young Russian Scientists" (grant MK-1378.2017.1).}}

\medskip

\medskip

\medskip

\medskip
\textbf{Nurlan Ismailov$^{a,b}$, Ivan Kaygorodov$^{c}$, Yury Volkov$^{d}$}
\medskip

{\tiny

$^{a}$ Universidade de S\~ao Paulo, IME, S\~ao Paulo, Brazil.

$^{b}$ Institute of Mathematics and Mathematical Modeling, Almaty, Kazakhstan.

$^{c}$ Universidade Federal do ABC, CMCC, Santo Andr\'{e}, Brazil.

$^{d}$ Saint Petersburg State University, Saint Petersburg, Russia.
\smallskip

    E-mail addresses:\smallskip
    
    Nurlan Ismailov (nurlan.ismail@gmail.com),
    
    Ivan Kaygorodov (kaygorodov.ivan@gmail.com),
    
    Yury Volkov (wolf86\_666@list.ru).

}

       \vspace{0.3cm}

\

{\bf Abstract.} 
We describe all degenerations of three dimensional anticommutative algebras $\mathfrak{Acom}_3$ and of three dimensional Leibniz algebras $\mathfrak{Leib}_3$ over $\mathbb{C}.$  
In particular, we describe all irreducible components and rigid algebras in the corresponding varieties.\smallskip

       \vspace{0.3cm}


\section{Introduction}

       \vspace{0.3cm}

Degenerations of algebras is an interesting subject, which was studied in various papers 
(see, for example, \cite{bu,BC99,S90,GRH,GRH2,BB09,BB14,GRH3,kpv16, kppv, kv16,  kv172, khud13,khom2, laur03,ra12,ale,ale2,aleis,casas13}).
In particular, there are many results concerning degenerations of algebras of low dimensions in a  variety defined by a set of identities.
One of important problems in this direction is the description of so-called rigid algebras. 
These algebras are of big interest, since the closures of their orbits under the action of the generalized linear group form irreducible components of a variety under consideration
(with respect to the Zariski topology). 
For example, rigid algebras were classified in the varieties of
low dimensional unital associative, Lie, Jordan and Leibniz algebras \cite{ikv17}.
There are fewer works in which the full information about degenerations was found for some variety of algebras.
This problem was solved 
for two dimensional pre-Lie algebras in \cite{BB09},  
for three dimensional Novikov algebras in \cite{BB14},  
for four dimensional Lie algebras in \cite{BC99}, 
for four dimensional Zinbiel algebras and nilpotent four-dimensional Leibniz algebras in \cite{kppv},
for nilpotent five  and six dimensional Lie algebras in \cite{S90,GRH}, 
for nilpotent five and six dimensional Malcev algebras in \cite{kpv16}, 
and for all two dimensional algebras in \cite{kv16}.


The most well known generalizations of Lie algebras are Leibniz, Malcev and binary Lie algebras.
The Leibniz algebras were introduced as a non-anticommutative generalization of Lie algebras.
The study of the structure theory and other properties of Leibniz algebras was initiated by Loday in \cite{lodaypir}.
An algebra $A$ is called a {\it  Leibniz  algebra}  if it satisfies the identity
$$(xy)z=(xz)y+x(yz).$$ 
The classification of all three dimensional Leibniz algebras can be found in \cite{ra12}.
Malcev algebras and binary Lie algebras are anticommutative. 
Gainov proved that there are no Malcev and binary Lie three dimensional algebras except Lie algebras \cite{G63}.
The description of all three dimensional anticommutative algebras was given in \cite{jap17}
and the central extensions of three dimensional anticommutative algebras were described in \cite{cfk18}.
In this paper we consider anticommutative algebras as a generalization of Lie algebras.
Note that some steps towards a classification of all three dimensional algebras have been done in \cite{Bekbaev}.

In this paper we give the full information about degenerations of three dimensional anticommutative and Leibniz algebras  over $\mathbb{C}.$
The vertices of this graph are the isomorphism classes of algebras in the variety under consideration.
An algebra $A$ degenerates to an algebra $B$ if and only if there is a path from the vertex corresponding to $A$ to the vertex corresponding to $B$.
We also describe rigid algebras and irreducible components in the corresponding varieties.

\section{Definitions and notation}

All spaces in this paper are considered over $\mathbb{C}$, and we write simply $dim$, $Hom$ and $\otimes$ instead of $dim_{\mathbb{C}}$, $Hom_{\mathbb{C}}$ and $\otimes_{\mathbb{C}}$. An algebra $A$ is a set with a structure of a vector space and a binary operation that induces a bilinear map from $A\times A$ to $A$.

Given an $n$ dimensional vector space $V$, the set $Hom(V \otimes V,V) \cong V^* \otimes V^* \otimes V$ 
is a vector space of dimension $n^3$. This space has a structure of the affine variety $\mathbb{C}^{n^3}.$ Indeed, let us fix a basis $e_1,\dots,e_n$ of $V$. Then any $\mu\in Hom(V \otimes V,V)$ is determined by $n^3$ structure constants $c_{i,j}^k\in\mathbb{C}$ such that
$\mu(e_i\otimes e_j)=\sum\limits_{k=1}^nc_{i,j}^ke_k$. A subset of $Hom(V \otimes V,V)$ is {\it Zariski-closed} if it can be defined by a set of polynomial equations in the variables $c_{i,j}^k$ ($1\le i,j,k\le n$).

Let $T$ be a set of polynomial identities.
All algebra structures on $V$ satisfying polynomial identities from $T$ form a Zariski-closed subset of the variety $Hom(V \otimes V,V)$. We denote this subset by $\mathbb{L}(T)$.
The general linear group $GL(V)$ acts on $\mathbb{L}(T)$ by conjugations:
$$ (g * \mu )(x\otimes y) = g\mu(g^{-1}x\otimes g^{-1}y)$$ 
for $x,y\in V$, $\mu\in \mathbb{L}(T)\subset Hom(V \otimes V,V)$ and $g\in GL(V)$.
Thus, $\mathbb{L}(T)$ is decomposed into $GL(V)$-orbits that correspond to the isomorphism classes of algebras. 
Let $O(\mu)$ denote the orbit of $\mu\in\mathbb{L}(T)$ under the action of $GL(V)$ and $\overline{O(\mu)}$ denote the Zariski closure of $O(\mu)$.

Let $A$ and $B$ be two $n$ dimensional algebras satisfying identities from $T$ and $\mu,\lambda \in \mathbb{L}(T)$ represent $A$ and $B$ respectively.
We say that $A$ degenerates to $B$ and write $A\to B$ if $\lambda\in\overline{O(\mu)}$.
Note that in this case we have $\overline{O(\lambda)}\subset\overline{O(\mu)}$. Hence, the definition of a degeneration does not depend on the choice of $\mu$ and $\lambda$. If $A\not\cong B$, then the assertion $A\to B$ is called a {\it proper degeneration}. We write $A\not\to B$ if $\lambda\not\in\overline{O(\mu)}$.

Let $A$ be represented by $\mu\in\mathbb{L}(T)$. Then  $A$ is  {\it rigid} in $\mathbb{L}(T)$ if $O(\mu)$ is an open subset of $\mathbb{L}(T)$.
 Recall that a subset of a variety is called irreducible if it cannot be represented as a union of two non-trivial closed subsets. A maximal irreducible closed subset of a variety is called an {\it irreducible component}.
It is well known that any affine variety can be represented as a finite union of its irreducible components in a unique way.
The algebra $A$ is rigid in $\mathbb{L}(T)$ if and only if $\overline{O(\mu)}$ is an irreducible component of $\mathbb{L}(T)$. 

We denote by $\mathfrak{ACom}_n$ the variety of $n$ dimensional anticommutative algebras and
by $\mathfrak{Leib}_n$ the variety of $n$ dimensional Leibniz algebras. 

We use the following notation: 

\begin{enumerate}
\item $Ann_L(A)=\{ a \in A \mid xa =0 \mbox{ for all } x\in A \}$ is the left  annihilator of $A;$
\item $A^{(+2)}$ is the space $\{xy+yx\mid x,y\in A\}$.
\end{enumerate}
Given spaces $U$ and $W$, we write simply $U>W$ instead of $dim\,U>dim\,W$.



\section{Methods} 

In the present work we use the methods applied to Lie algebras in \cite{BC99,GRH,GRH2,S90}.
First of all, if $A\to B$ and $A\not\cong B$, then $Der(A)<Der(B)$, where $Der(A)$ is the Lie algebra of derivations of $A$. We will compute the dimensions of algebras of derivations and will check the assertion $A\to B$ only for such $A$ and $B$ that $Der(A)<Der(B)$. Secondly, if $A\to C$ and $C\to B$ then $A\to B$. If there is no $C$ such that $A\to C$ and $C\to B$ are proper degenerations, then the assertion $A\to B$ is called a {\it primary degeneration}. If $Der(A)<Der(B)$ and there are no $C$ and $D$ such that $C\to A$, $B\to D$, $C\not\to D$ and one of the assertions $C\to A$ and $B\to D$ is a proper degeneration,  then the assertion $A \not\to B$ is called a {\it primary non-degeneration}. It suffices to prove only primary degenerations and non-degenerations to describe degenerations in the variety under consideration. It is easy to see that any algebra degenerates to the algebra with zero multiplication. From now on we use this fact without mentioning it.


To prove primary degenerations, we will construct families of matrices parametrized by $t$. Namely, let $A$ and $B$ be two algebras represented by the structures $\mu$ and $\lambda$ from $\mathbb{L}(T)$ respectively. Let $e_1,\dots, e_n$ be a basis of $V$ and $c_{i,j}^k$ ($1\le i,j,k\le n$) be the structure constants of $\lambda$ in this basis. If there exist $a_i^j(t)\in\mathbb{C}$ ($1\le i,j\le n$, $t\in\mathbb{C}^*$) such that $E_i^t=\sum\limits_{j=1}^na_i^j(t)e_j$ ($1\le i\le n$) form a basis of $V$ for any $t\in\mathbb{C}^*$, and the structure constants of $\mu$ in the basis $E_1^t,\dots, E_n^t$ are such polynomials $c_{i,j}^k(t)\in\mathbb{C}[t]$ that $c_{i,j}^k(0)=c_{i,j}^k$, then $A\to B$. In this case  $E_1^t,\dots, E_n^t$ is called a {\it parametrized basis} for $A\to B$.

Note also the following fact. Let $B(\alpha)$ be a series of algebras parametrized by $\alpha\in\mathbb{C}$ and $e_1,\dots,e_n$ be a basis of $V$. Suppose also that, for any $\alpha\in\mathbb{C}$, the algebra $B(\alpha)$ can be represented by a structure $\mu(\alpha)\in\mathbb{L}(T)$ having structure constants $c_{i,j}^k(\alpha)\in\mathbb{C}$ in the basis $e_1,\dots,e_n$, where $c_{i,j}^k(t)\in\mathbb{C}[t]$ for all $1\le i,j,k\le n$. Let $A$ be an algebra such that $A\to B(\alpha)$ for $\alpha\in\mathbb{C}\setminus S$, where $S$ is a finite subset of $\mathbb{C}$. Then $A\to B(\alpha)$ for all $\alpha\in\mathbb{C}$. Indeed, if $\lambda\in \mathbb{L}(T)$ represents $A$, then we have $\mu(\alpha)\in\overline{\{\mu(\beta)\}_{\beta\in\mathbb{C}\setminus S}}\subset \overline{O(\lambda)}$ for any $\alpha\in\mathbb{C}$. Thus, to prove that $A\to B(\alpha)$ for all $\alpha\in\mathbb{C}$ we will construct degenerations that are valid for all but finitely many $\alpha$.

Let us describe the methods for proving primary non-degenerations. The main tool for this is the following lemma.

\begin{Lem}[\cite{GRH}]\label{main}
Let $\mathcal{B}$ be a Borel subgroup of $GL(V)$ and $\mathcal{R}\subset \mathbb{L}(T)$ be a $\mathcal{B}$-stable closed subset.
If $A \to B$ and $A$ can be represented by $\mu\in\mathcal{R}$ then there is $\lambda\in \mathcal{R}$ that represents $B$.
\end{Lem}

In particular, it follows from Lemma \ref{main} that $A\not\to B$ in the following cases:
\begin{enumerate}
\item $Ann_L(A)>Ann_L(B)$;
\item $A^{(+2)}<B^{(+2)}$.
\end{enumerate}
In the cases where all of these criteria cannot be applied to prove $A\not\to B$, we will define $\mathcal{R}$ by a set of polynomial equations and will give a basis of $V$, in which the structure constants of $\mu$ give a solution to all these equations. We will omit everywhere the verification of the fact that $\mathcal{R}$ is stable under the action of the subgroup of upper triangular matrices and of the fact that $\lambda\not\in\mathcal{R}$ for any choice of a basis of $V$. These verifications can be done by direct calculations.

If the number of orbits under the action of $GL(V)$ on  $\mathbb{L}(T)$ is finite, then the graph of primary degenerations gives the whole picture. In particular, the description of rigid algebras and irreducible components can be easily obtained. Since the variety $\mathfrak{Leib}_3$ contains infinitely many non-isomorphic algebras, we have to fulfill some additional work. Let $A(*):=\{A(\alpha)\}_{\alpha\in I}$ be a set of algebras, and let $B$ be another algebra. Suppose that, for $\alpha\in I$, $A(\alpha)$ is represented by the structure $\mu(\alpha)\in\mathbb{L}(T)$ and $B\in\mathbb{L}(T)$ is represented by the structure $\lambda$. Then $A(*)\to B$ means $\lambda\in\overline{\{O(\mu(\alpha))\}_{\alpha\in I}}$, and $A(*)\not\to B$ means $\lambda\not\in\overline{\{O(\mu(\alpha))\}_{\alpha\in I}}$.

Let $A(*)$, $B$, $\mu(\alpha)$ ($\alpha\in I$) and $\lambda$ be as above. To prove $A(*)\to B$ it is enough to construct a family of pairs $(f(t), g(t))$ parametrized by $t\in\mathbb{C}^*$, where $f(t)\in I$ and $g(t)\in GL(V)$. Namely, let $e_1,\dots, e_n$ be a basis of $V$ and $c_{i,j}^k$ ($1\le i,j,k\le n$) be the structure constants of $\lambda$ in this basis. If we construct $a_i^j:\mathbb{C}^*\to \mathbb{C}$ ($1\le i,j\le n$) and $f: \mathbb{C}^* \to I$ such that $E_i^t=\sum\limits_{j=1}^na_i^j(t)e_j$ ($1\le i\le n$) form a basis of $V$ for any  $t\in\mathbb{C}^*$, and the structure constants of $\mu_{f(t)}$ in the basis $E_1^t,\dots, E_n^t$ are such polynomials $c_{i,j}^k(t)\in\mathbb{C}[t]$ that $c_{i,j}^k(0)=c_{i,j}^k$, then $A(*)\to B$. In this case  $E_1^t,\dots, E_n^t$ and $f(t)$ are called a parametrized basis and a {\it parametrized index} for $A(*)\to B$ respectively.

We now explain how to prove $A(*)\not\to B$. Note that if $dim\,Der(A(\alpha))>dim\,Der(B)$ for all $\alpha\in I$ then $A(*)\not\to B$. One can use also the following generalization of Lemma \ref{main}, whose proof is the same as the proof of Lemma \ref{main}.

\begin{Lem}\label{gmain}
Let $\mathcal{B}$ be a Borel subgroup of $GL(V)$ and $\mathcal{R}\subset \mathbb{L}(T)$ be a $\mathcal{B}$-stable closed subset.
If $A(*) \to B$ and for any $\alpha\in I$ the algebra $A(\alpha)$ can be represented by a structure $\mu(\alpha)\in\mathcal{R}$, then there is $\lambda\in \mathcal{R}$ representing $B$.
\end{Lem}

\section{Classification and degenerations of three dimensional anticommutative algebras}

First we consider the variety $\mathfrak{ACom}_3$. Let us fix the basis $e_1,e_2,e_3$ of $V$. Any structure $\mu\in\mathfrak{ACom}_3$ with structure constants $c_{i,j}^k$ ($1\le i,j,k\le 3$) is determined by the $3\times 3$ matrix $A^{\mu}$ whose $(i,j)$-entry is $(-1)^{i-1}c_{u,v}^j$, where $(u,v)$ is a unique pair of numbers such that $u,v\in\{1,2,3\}\setminus \{i\}$ and $u<v$. Since $\mathbb{C}$ is an algebraically closed field the structure $\lambda$ belongs to $O(\mu)$ if and only if there is a nonsingular matrix $X$ such that $A^{\lambda}=X^TA^{\mu}X$ by \cite[Corollary 2.4]{jap17}. Then the classification of three dimensional anticommutative algebras modulo isomorphism can be obtained from the classification of bilinear forms modulo congruence given in \cite{ST}.

We denote by $W$ some four dimensional space that contains $V$ as a subspace and by $e_4$ some fixed vector of $W$ such that $W=V\oplus\mathbb{C}e_4$.
Let us now consider four dimensional algebras $A$ such that $A(A^2)=(A^2)A=0$ and $dim\,A^2\le 1$. It is easy to see that such an algebra can be represented by a structure $\chi$ on $W$ such that $\chi(W,W)\subset\mathbb{C}e_4$ and $\chi(W,e_4)=\chi(e_4,W)=0$.
Such a structure is defined by the $3\times 3$ matrix $B^{\chi}$, whose $(i,j)$-entry is $d_{i,j}^4$, where $d_{i,j}^k$ ($1\le i,j,k\le 4$) are the structure constants of $\chi$. It is clear that two such structures $\chi$ and $\eta$ lie in the same orbit if and only if there is a nonsingular matrix $X$ such that $B^{\eta}=X^TB^{\chi}X$.
Now we put in correspondence to an anticommutative algebra structure $\mu$ on $V$ the structure $\chi_{\mu}$ on $W$ satisfying the properties above with $B^{\chi_{\mu}}=A^{\mu}$. As it was explained above, we get a bijection between orbits of $\mathfrak{ACom}_3$ and isomorphism classes of four dimensional algebras $A$ such that $A(A^2)=(A^2)A=0$ and $dim\,A^2\le 1$. Moreover, it is clear that if $\lambda\in \overline{O(\mu)}$, then $\chi_{\lambda}\in \overline{O(\chi_{\mu})}$. The converse assertion follows from \cite[Proposition 1.7]{GRH} and the fact that the set of structures $\chi$ on $W$ satisfying $\chi(W,W)\subset\mathbb{C}e_4$ and $\chi(W,e_4)=\chi(e_4,W)=0$ is a closed subset stable under the action of lower triangular matrices.

Thus, isomorphism classes and degenerations of three dimensional anticommutative algebras can be transfered from the isomorphism classes and degenerations of four dimensional algebras $A$ such that $A(A^2)=(A^2)A=0$ and $dim\,A^2\le 1$. The last mentioned problem is a part of the problems that were solved in \cite{fipe1, kppv}. Unfortunately both of the mentioned works have problems. Some degenerations are missed in the paper \cite{fipe1}. All degenerations between algebras that we are interested in are described correctly in \cite{kppv}, but the classification used in this paper lost one algebra.
In the current work we will use the results of \cite{kppv}. Let us first deduce the classification of three dimensional anticommutative algebras using the last mentioned paper.
We do this in Table A.1 below, where in the first column we put the names of anticommutative three dimensional algebras, in the second column we put the corresponding names of algebras from \cite{kppv}, in the third and fourth columns we put multiplication tables and dimensions of algebras of derivations of the corresponding anticommutative algebras. We omit products of basic elements whose values are zero or can be recovered from the anticommutativity and given products.
Note that if $dim\,Der(\mu)=k$, then $dim\,Der(\chi_{\mu})=k+4$.

Here $\mathfrak{g}_{1}$, $\mathfrak{g}_{2}$, $\mathfrak{g}_{3}^{\alpha}$ and $\mathfrak{g}_4$ are Lie algebras and $\mathcal{A}_2$ corresponds to the algebra missed in \cite{kppv} that is denoted by $\mathfrak{N}_0$ in this paper. We have $\mathfrak{g}^{\alpha}_3\cong \mathfrak{g}^{\beta}_3$ and $\mathcal{A}^{\alpha}_1\cong \mathcal{A}^{\beta}_1$ if $\alpha\beta=1$ and there are no other nontrivial isomorphisms between the algebras in the table. All degenerations and non-degenerations between the algebras from the column $\mathbb{B}$ that do not involve $\mathfrak{N}_0$ are described in \cite{kppv}. Thus, it remains to describe degenerations involving $\mathcal{A}_2$.

Note that $dim\,Ann_L(\mathfrak{N}_0)>dim\,Ann_L(\mathfrak{N}_3(0))=dim\,Ann_L(\mathfrak{N}_{10})$ and $dim\,\mathfrak{Der}(\mathfrak{N}_0)<dim\,\mathfrak{Der}(\mathfrak{N}_3(0))=dim\,\mathfrak{Der}(\mathfrak{N}_{10})$, and hence there are no degenerations between $\mathfrak{N}_0$, $\mathfrak{N}_3(0)$ and $\mathfrak{N}_{10}$. Note now that, for any $\alpha\in\mathbb{C}$ we have
\begin{itemize}
\item the degeneration $\mathcal{A}_2\to \mathfrak{g}^{\alpha}_3$ given by the parametrized basis 
$$E_1^t=te_3, E_2^t=te_1, E_3^t=e_1+(\alpha+t)e_2+e_3;$$
\item the degeneration $\mathcal{A}_1^{\alpha}\to\mathcal{A}_2$ given by the parametrized basis
$$E_1^t=te_2, E_2^t=-e_1, E_3^t=\alpha e_1-e_2+e_3.$$
\end{itemize}

Thus, we get the following result.

\begin{Th}
\label{acom3}
The graph of primary degenerations for  $\mathfrak{ACom}_3$ can be obtained from the graph given in Figure 1 below by deleting all vertices named with $\mathcal{L}$.
\end{Th}

Since $\mathfrak{ACom}_3$ is isomorphic to $\mathbb{C}^{9}$ as an algebraic variety, it is irreducible and equals $\overline{ O(\mathcal{A}_1^{\alpha} )  }$

\section{Degenerations of three dimensional Leibniz algebras}
The classification of three dimensional non-Lie Leibniz algebras is presented in Table A.2 below. 


\normalsize 

\begin{Th}\label{degLei}
\label{leib3}
The graph of primary degenerations for  $\mathfrak{Leib}_3$ can be obtained from the graph given in Figure 1 below by deleting all vertices named with $\mathcal{A}$.
\end{Th}

\begin{Proof}  
We prove all the required primary degenerations in Table A.3 below. 
Let us consider the degeneration $\mathfrak{L}^{\beta}_1  \to  \mathfrak{L}_{2}$ to clarify our formulas. 
Write nonzero products in $\mathfrak{L}^{\beta}_1 $ in the basis  
$E_i^t$: 
$$E_2^tE_2^t=\beta t^2E_1^t, \  E_3^tE_2^t= tE_1^t, E_3^tE_3^t=E_1^t.$$
It is easy to see that for $t=0$ we obtain the multiplication table of $\mathfrak{L}_{2}$. The remaining degenerations can be interpreted in the same way.

A part of non-degenerations is given in Table A.4 below. Whenever an algebra named by $A$ with the basis $f_1,f_2,f_3$ appear in this table, $c_{i,j}^k$ ($1\le i,j,k\le 3$) denote the structure constants of $A$ in the given basis and $A_i$ ($1\le i\le 3$) denotes the subspace of $A$ generated by $f_i,\ldots,f_3$.

In the rest of the proof we will use ideas from the proof of \cite[Theorem 2]{ikv17}.
All the remaining degenerations involve only solvable non-nilpotent Leibniz algebras with a two dimensional nilpotent radical. 
Moreover,  each of them is represented in Table A.2 by a structure 
$\mu$ such that $\langle e_1,e_2\rangle$ is the nilpotent radical and the structure constants $c_{ij}^k$ ($1\le i,j,k\le 3$) satisfy the conditions $c_{ij}^k=0$ if $i,j\le 2$ and $k\ge \min(i,j)$ and $c_{3i}^j=c_{i3}^j=0$ for any $1\le i,j\le 2$ such that $j<i$. During this proof we will call a structure with three dimensional nilpotent radical that satisfies the described conditions a {\it standard structure}.
Let us put in correspondence to a standard structure $\mu$ the $4$-tuple $S_{\mu}=(c_{13}^1,c_{31}^1,c_{23}^2,c_{32}^2)\in\mathbb{C}^4$. It is not difficult to show that if $S_{\mu}=(a_1,b_1,a_2,b_2)$ and $\lambda\in O(\mu)$ is a standard structure, then there is some permutation $\sigma:\{1,2\}\rightarrow\{1,2\}$ and some $c\in\mathbb{C}^*$ such that $S_{\lambda}=(ca_{\sigma(1)},cb_{\sigma(1)},ca_{\sigma(2)},cb_{\sigma(2)})$. 
Suppose now that $\{\mu_s\}_{s\in T}$ is a set of standard structures, $S_{\mu_s}=(a_{1,s},b_{1,s},a_{2,s},b_{2,s})$, and the homogeneous linear polynomials $f_1,\dots,f_l\in\mathbb{C}[x_1,x_2,x_3,x_4]$ are such that $f_r(a_{1,s},b_{1,s},a_{2,s},b_{2,s})=0$ for all $s\in T$ and $1\le r\le l$. If $\lambda$ is a standard structure with $S_{\lambda}=(a_1,b_1,a_2,b_2)$, then it easily follows from Lemma \ref{gmain} that there is some permutation $\sigma:\{1,2\}\rightarrow\{1,2\}$ and some $c\in\mathbb{C}^*$ such that $f_r(ca_{\sigma(1)},cb_{\sigma(1)},ca_{\sigma(2)},cb_{\sigma(2)})=0$ for all $1\le r\le l$.
Thus, we get

$$
\mathfrak{B}\not\to\mathfrak{L}^0_6  \hspace{0.2cm} \mbox{for} \hspace{0.2cm}   \mathfrak{B}\in\{\mathfrak{L}_6^{\alpha\neq0},\mathfrak{L}_{7}\}\hspace{0.8cm}\mbox{ and }\hspace{0.8cm}
\mathfrak{B}\not\to\mathfrak{L}^1_6  \hspace{0.2cm} \mbox{for} \hspace{0.2cm}   \mathfrak{B}\in\{\mathfrak{L}_4^{\alpha},\mathfrak{L}_5,\mathfrak{L}_6^{\alpha\neq1},\mathfrak{L}_{9}\}.
$$
\end{Proof}

\begin{corollary}\label{ir_leib3} The irreducible components  of $\mathfrak{Leib}_3$ are
$$
\begin{aligned}
\mathcal{C}_1   &=\overline{ O(\{\mathfrak{g}_3^{\alpha}\}_{\alpha\in\mathbb{C}} )  }=  \{ \mathfrak{g}_1, \mathfrak{g}_2,  \mathfrak{g}_3^{\alpha},  \mathbb{C}^3 \}_{\alpha\in\mathbb{C}};\\
\mathcal{C}_2   &=\overline{ O(\mathfrak{g}_4 )  }=  \{ \mathfrak{g}_1,
\mathfrak{g}_3^{-1}, \mathfrak{g}_4, \mathbb{C}^3 \};\\
\mathcal{C}_3   &=\overline{ O(\{\mathfrak{L}_4^{\alpha}\}_{\alpha\in\mathbb{C}} )  }=  
\{ \mathfrak{g}_1, \mathfrak{g}_3^0, \mathfrak{L}_1^{\beta}, \mathfrak{L}_2,  \mathfrak{L}_3, \mathfrak{L}_4^{\alpha},  \mathfrak{L}_6^0,  \mathbb{C}^3 \}_{\alpha,\beta\in\mathbb{C}};\\
\mathcal{C}_4   &=\overline{ O(\mathfrak{L}_5 )  }=  
\{ \mathfrak{L}_2, \mathfrak{L}_3, \mathfrak{L}_4^2, \mathfrak{L}_5,   \mathbb{C}^3 \};\\
\mathcal{C}_5   &=\overline{ O(\{\mathfrak{L}_6^{\alpha}\}_{\alpha\in\mathbb{C}} )  }=  
\{\mathfrak{L}_1^0, \mathfrak{L}_2, \mathfrak{L}_6^{\alpha},   \mathfrak{L}_7, \mathfrak{L}_{8}, \mathfrak{L}_{9}, \mathbb{C}^3 \}_{\alpha\in\mathbb{C}}.\\
\end{aligned}
$$ 
In particular, the set of rigid algebras in the variety $\mathfrak{Leib}_3$ is formed by $\mathfrak{g}_4$ and $\mathfrak{L}_5$.
\end{corollary} 
\begin{Proof} All degenerations and non-degenerations that do not follow directly from Theorem \ref{degLei} follow from Talble A.5.
\end{Proof}

\section*{Figure 1: 
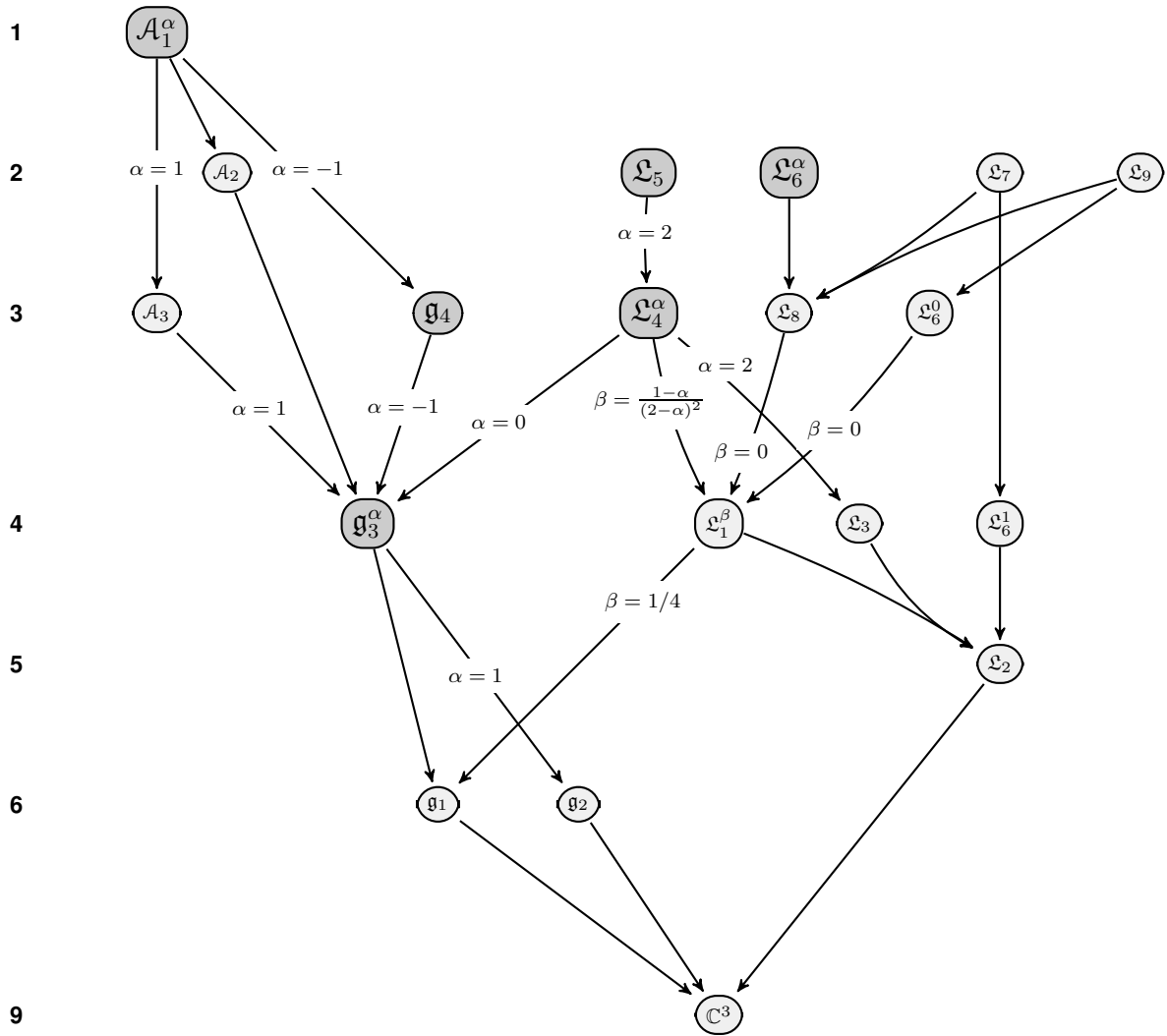
The graph of primary degenerations for \\  Lie, anticommutative and Leibniz three dimensional algebras.
}

\begin{center}

\

\

\begin{tikzpicture}[->,>=stealth',shorten >=0.05cm, auto, node distance=0.95cm, thick,
                    main  node/.style={rectangle,draw,fill=gray!12, rounded corners=1.5ex,font=\sffamily \tiny \bfseries },
                    rigid node/.style={rectangle,draw,fill=black!20,rounded corners=1.5ex,font=\sffamily \bf \bfseries }, 
                    ver  node/.style={rectangle,draw,fill=gray!12, rounded corners=1.5ex,font=\sffamily \bf \bfseries },
                    ser node/.style={rectangle,draw,fill=green!20,rounded corners=1.5ex,font=\sffamily \bf \bfseries },
                    nrig node/.style={rectangle,draw,fill=cyan!20,rounded corners=1.5ex,font=\sffamily \tiny \bfseries },
                    nrigdead node/.style={rectangle,draw,fill=teal!60,rounded corners=1.5ex,font=\sffamily \tiny \bfseries }, 
                    lie node/.style={rectangle,draw,fill=yellow!20,rounded corners=1.5ex,font=\sffamily \tiny \bfseries },
                    lier node/.style={rectangle,draw,fill=yellow!80,rounded corners=1.5ex,font=\sffamily \tiny \bfseries },
                    style={draw,font=\sffamily \scriptsize \bfseries }]
\node (0)   {};

\node (00a1) [right of=0] {};
\node (00a2) [right of=00a1] {};
\node (00a3) [right of=00a2] {};
\node (00a4) [right of=00a3] {};
\node (00a5) [right of=00a4] {};
\node (00a6) [right of=00a5] {};
\node (00a7) [right of=00a6] {};
\node (00a8) [right of=00a7] {};
\node (00a9) [right of=00a8] {};
\node (00a10) [right of=00a9] {};
\node (00a11) [right of=00a10] {};
\node (00a12) [right of=00a11] {};
\node (00a13) [right of=00a12] {};
\node (00a14) [right of=00a13] {};
\node (00a15) [right of=00a14] {};
\node (00a16) [right of=00a15] {};
\node (00a17) [right of=00a16] {};
\node (00a18) [right of=00a17] {};
\node (00a19) [right of=00a18] {};
\node (00a20) [right of=00a19] {};
\node (00a21) [right of=00a20] {};
\node (00a22) [right of=00a21] {};
\node (00a23) [right of=00a22] {};
\node (00a24) [right of=00a23] {};
\node (00a25) [right of=00a24] {};
\node (00a26) [right of=00a25] {};
\node (00a27) [right of=00a26] {};

\node (01) [below of=0] {1};

\node (01a1) [right of=01] {};
\node (01a2) [right of=01a1] {};
\node (01a3) [right of=01a2] {};
\node (01a4) [right of=01a3] {};
\node (01a5) [right of=01a4] {};
\node (01a6) [right of=01a5] {};
\node (01a7) [right of=01a6] {};
\node (01a8) [right of=01a7] {};
\node (01a9) [right of=01a8] {};
\node (01a10) [right of=01a9] {};
\node (01a11) [right of=01a10] {};
\node (01a12) [right of=01a11] {};
\node (01a13) [right of=01a12] {};
\node (01a14) [right of=01a13] {};
\node (01a15) [right of=01a14] {};
\node (01a16) [right of=01a15] {};
\node (01a17) [right of=01a16] {};

\node (021)[below of=01]{};
\node (02) [below of=021] {2};

\node (02a1) [right of=02] {};
\node (02a2) [right of=02a1] {};
\node (02a3) [right of=02a2] {};
\node (02a4) [right of=02a3] {};
\node (02a5) [right of=02a4] {};
\node (02a6) [right of=02a5] {};
\node (02a7) [right of=02a6] {};
\node (02a8) [right of=02a7] {};
\node (02a9) [right of=02a8] {};
\node (02a10) [right of=02a9] {};
\node (02a11) [right of=02a10] {};
\node (02a12) [right of=02a11] {};
\node (02a13) [right of=02a12] {};
\node (02a14) [right of=02a13] {};
\node (02a15) [right of=02a14] {};
\node (02a16) [right of=02a15] {};
\node (02a17) [right of=02a16] {};
\node (02a18) [right of=02a17] {};
\node (02a19) [right of=02a18] {};
\node (02a20) [right of=02a19] {};
\node (02a21) [right of=02a20] {};
\node (02a22) [right of=02a21] {};
\node (02a23) [right of=02a22] {};
\node (02a24) [right of=02a23] {};
\node (02a25) [right of=02a24] {};
\node (02a26) [right of=02a25] {};
\node (02a27) [right of=02a26] {};

\node (031)[below of=02]{};
\node (03)[below of=031]{3};

\node (03a1) [right of=03] {};
\node (03a2) [right of=03a1] {};
\node (03a3) [right of=03a2] {};
\node (03a4) [right of=03a3] {};
\node (03a5) [right of=03a4] {};
\node (03a6) [right of=03a5] {};
\node (03a7) [right of=03a6] {};
\node (03a8) [right of=03a7] {};
\node (03a9) [right of=03a8] {};
\node (03a10) [right of=03a9] {};
\node (03a11) [right of=03a10] {};
\node (03a12) [right of=03a11] {};
\node (03a13) [right of=03a12] {};
\node (03a14) [right of=03a13] {};
\node (03a15) [right of=03a14] {};
\node (03a16) [right of=03a15] {};
\node (03a17) [right of=03a16] {};
\node (03a18) [right of=03a17] {};
\node (03a19) [right of=03a18] {};
\node (03a20) [right of=03a19] {};
\node (03a21) [right of=03a20] {};
\node (03a22) [right of=03a21] {};
\node (03a23) [right of=03a22] {};
\node (03a24) [right of=03a23] {};
\node (03a25) [right of=03a24] {};
\node (03a26) [right of=03a25] {};
\node (03a27) [right of=03a26] {};
*\node (03a28) [right of=03a27] {};

\node (041)[below of=03]{};
\node (042)[below of=041]{};

\node (04) [below of=042] {4};

\node (04a1) [right of=04] {};
\node (04a2) [right of=04a1] {};
\node (04a3) [right of=04a2] {};
\node (04a4) [right of=04a3] {};
\node (04a5) [right of=04a4] {};
\node (04a6) [right of=04a5] {};+
\node (04a7) [right of=04a6] {};
\node (04a8) [right of=04a7] {};
\node (04a9) [right of=04a8] {};
\node (04a10) [right of=04a9] {};
\node (04a11) [right of=04a10] {};
\node (04a12) [right of=04a11] {};
\node (04a13) [right of=04a12] {};
\node (04a14) [right of=04a13] {};
\node (04a15) [right of=04a14] {};
\node (04a16) [right of=04a15] {};
\node (04a17) [right of=04a16] {};
\node (04a18) [right of=04a17] {};
\node (04a19) [right of=04a18] {};
\node (04a20) [right of=04a19] {};
\node (04a21) [right of=04a20] {};
\node (04a22) [right of=04a21] {};
\node (04a23) [right of=04a22] {};
\node (04a24) [right of=04a23] {};
\node (04a25) [right of=04a24] {};
\node (04a26) [right of=04a25] {};
\node (04a27) [right of=04a26] {};
\node (04a28) [right of=04a27] {};

\node (051)[below of=04]{};
\node (05) [below of=051] {5};

\node (05a1) [right of=05] {};
\node (05a2) [right of=05a1] {};
\node (05a3) [right of=05a2] {};
\node (05a4) [right of=05a3] {};
\node (05a5) [right of=05a4] {};
\node (05a6) [right of=05a5] {};
\node (05a7) [right of=05a6] {};
\node (05a8) [right of=05a7] {};
\node (05a9) [right of=05a8] {};
\node (05a10) [right of=05a9] {};
\node (05a11) [right of=05a10] {};
\node (05a12) [right of=05a11] {};
\node (05a13) [right of=05a12] {};
\node (05a14) [right of=05a13] {};
\node (05a15) [right of=05a14] {};
\node (05a16) [right of=05a15] {};
\node (05a17) [right of=05a16] {};
\node (05a18) [right of=05a17] {};
\node (05a19) [right of=05a18] {};
\node (05a20) [right of=05a19] {};
\node (05a21) [right of=05a20] {};
\node (05a22) [right of=05a21] {};
\node (05a23) [right of=05a22] {};
\node (05a24) [right of=05a23] {};
\node (05a25) [right of=05a24] {};
\node (05a26) [right of=05a25] {};
\node (05a27) [right of=05a26] {};
\node (05a28) [right of=05a27] {};

\node (061)[below of=05]{};
\node (06)[below of=061] {6};

\node (06a1) [right of=06] {};
\node (06a2) [right of=06a1] {};
\node (06a3) [right of=06a2] {};
\node (06a4) [right of=06a3] {};
\node (06a5) [right of=06a4] {};
\node (06a6) [right of=06a5] {};
\node (06a7) [right of=06a6] {};
\node (06a8) [right of=06a7] {};
\node (06a9) [right of=06a8] {};
\node (06a10) [right of=06a9] {};
\node (06a11) [right of=06a10] {};
\node (06a12) [right of=06a11] {};
\node (06a13) [right of=06a12] {};
\node (06a14) [right of=06a13] {};
\node (06a15) [right of=06a14] {};
\node (06a16) [right of=06a15] {};
\node (06a17) [right of=06a16] {};
\node (06a18) [right of=06a17] {};
\node (06a19) [right of=06a18] {};
\node (06a20) [right of=06a19] {};
\node (06a21) [right of=06a20] {};
\node (06a22) [right of=06a21] {};
\node (06a23) [right of=06a22] {};
\node (06a24) [right of=06a23] {};
\node (06a25) [right of=06a24] {};
\node (06a26) [right of=06a25] {};
\node (06a27) [right of=06a26] {};
\node (06a28) [right of=06a27] {};

\node (071)[below of=06]{};
\node (072)[below of=071]{};
\node (07)[below of=072] {9};

\node (07a1) [right of=07] {};
\node (07a2) [right of=07] {};
\node (07a3) [right of=07a2] {};
\node (07a4) [right of=07a3] {};
\node (07a5) [right of=07a4] {};
\node (07a6) [right of=07a5] {};
\node (07a7) [right of=07a6] {};
\node (07a8) [right of=07a7] {};
\node (07a9) [right of=07a8] {};
\node (07a10) [right of=07a9] {};
\node (07a11) [right of=07a10] {};
\node (07a12) [right of=07a11] {};
\node (07a13) [right of=07a12] {};
\node (07a14) [right of=07a13] {};
\node (07a15) [right of=07a14] {};
\node (07a16) [right of=07a15] {};
\node (07a17) [right of=07a16] {};
\node (07a18) [right of=07a17] {};
\node (07a19) [right of=07a18] {};
\node (07a20) [right of=07a19] {};
\node (07a21) [right of=07a20] {};
\node (07a22) [right of=07a21] {};
\node (07a23) [right of=07a22] {};
\node (07a24) [right of=07a23] {};
\node (07a25) [right of=07a24] {};
\node (07a26) [right of=07a25] {};
\node (07a27) [right of=07a26] {};
\node (07a28) [right of=07a27] {};

\node   (g5l)     [right of =03a4] {};
\node   [rigid node] (g4)     [right of =g5l] {$\mathfrak{g}_{4}$};
\node   [rigid node] (g31)     [right of =04a4] {$\mathfrak{g}^{\alpha}_{3}$};
\node   [main node] (g1)     [right of =06a5] {$\mathfrak{g}_{1}$};
\node   [main node] (g32)     [right of =06a7] {$\mathfrak{g}_{2}$};

\node   [rigid node] (l5)       [right of =02a8] {$\mathfrak{L}_{5}$};
\node   [rigid node] (l61)      [right of =02a10] {$\mathfrak{L}^{\alpha}_{6}$};
\node   [main node] (l7)        [right of =02a13] {$\mathfrak{L}_{7}$};
\node   [main node] (l9)        [right of =02a15] {$\mathfrak{L}_{9}$};

\node   [rigid node] (l4)      [right of =03a8] {$\mathfrak{L}^{\alpha}_{4}$};
\node   [main node] (l8)      [right of =03a10] {$\mathfrak{L}_{8}$};
\node   [main node] (l60)     [right of =03a12] {$\mathfrak{L}^{0}_{6}$};

\node   [main node] (l62)     [right of =04a13] {$\mathfrak{L}^{1}_{6}$};
\node   [main node] (l1)     [right of =04a9] {$\mathfrak{L}^{\beta}_{1}$};
\node   [main node] (l3)      [right of =04a11] {$\mathfrak{L}_3$};

\node   [main node] (l2)     [right of =05a13] {$\mathfrak{L}_2$};

\node   [rigid node] (ant5)     [right of =01a1] {$\mathcal{A}_1^{\alpha}$};
\node   (ant6l)     [right of =02a1] {};
\node   [main node] (ant6)     [right of =ant6l] {$\mathcal{A}_2$};
\node   [main node] (ant8)     [right of =03a1] {$\mathcal{A}_3$};

\node  [main node] (t20) [right of=07a10] {$\mathbb{C}^3$};

\path[every node/.style={font=\sffamily\small}]

(g1) edge        node[above] {} (t20)
(g32) edge     node[above] {} (t20)
(g31) edge  node[above] {} (g1)
(g31) edge   node[above=-10, fill=white]{\tiny $\alpha=1$} (g32)
(g4) edge   node[above=-4, fill=white]{\tiny $\alpha=-1$} (g31)

(ant5) edge   node[above=-4, fill=white]{\tiny $\alpha=1$} (ant8)
(ant5) edge   node{} (ant6)
(ant5) edge  node[above=2, right=-15, fill=white]{\tiny $\alpha=-1$} (g4)

(ant6) edge [bend left=0]  node{} (g31)
(ant8) edge   node[above=-4, fill=white]{\tiny $\alpha=1$} (g31)

(l1) edge   [bend left=0]  node[above=25, right=7, fill=white]{\tiny $\beta=1/4$} (g1)
(l1) edge    [bend left=5] node[above] {} (l2)

(l2) edge   node[above] {} (t20)

(l3) edge    [bend left=-15] node[above] {} (l2)

(l4)  edge    [bend left=-8] node[above=7, right=-35, fill=white]{\tiny $\beta=\frac{1-\alpha}{(2-\alpha)^2}$}  (l1)

(l4) edge    [bend left=5] node[above=20,right=-30,fill=white]{\tiny $\alpha=2$} (l3)


(l4) edge    [bend left=0] node[above=-1, right=-18, fill=white]{{\tiny $\alpha=0$}} (g31)

(l5) edge   [bend right=4] node[above=-3, fill=white]{\tiny $\alpha=2$} (l4)

(l60) edge   [bend left=8] node[above=-1, right=-15,, fill=white]{{\tiny $\beta=0$}} 
(l1)


(l61)  edge  node[above]{}  (l8)
(l62) edge   node[above] {} (l2)

(l7) edge   [bend right=0] node{} (l62)
(l7) edge   [bend left=6] node[above] {} (l8)


(l8) edge [bend left=4]  node[above=-14, right=-22, fill=white]{{\tiny $\beta=0$}} (l1)

(l9) edge   [bend left=-6] node[above] {} (l8)
(l9) edge   node[above] {} (l60);

\end{tikzpicture}

\end{center}

\newpage 
{\bf Appendix A. Tables}
\begin{center}  \footnotesize 
$$
\begin{array}{|l|c|lrrrrl|l|c|} 

\multicolumn{10}{l}{ \mbox{{\bf Table A.1.} 
{\it Classification of three dimensional anticommutative algebras.}}}\\ 

\multicolumn{10}{l}{}\\
\hline

\mathbb{A}& 
\multicolumn{7}{|l|}{\mathbb{B}} &\mbox{ Multiplication tables } & \mathfrak{Der} \\ 
\hline

\mathfrak{g}_{1} & 
\multicolumn{7}{|l|}{ \begin{array}{l} \mathfrak{N}_1^{\mathbb{C}^2} \end{array} }& \begin{array}{c}e_2e_3 =e_1 \end{array} & 6\\ 
\hline

\mathfrak{g}_2 &
\multicolumn{7}{|l|}{\begin{array}{l} \mathfrak{N}_1^{\mathbb{C}} \end{array}}& \begin{array}{c}e_1e_3 =e_1, \  e_2e_3=e_2 \end{array}& 6\\ 
\hline

\mathfrak{g}^{\alpha}_3& \multicolumn{7}{|l|}{
\begin{array}{l}\mathfrak{N}_3^{\mathbb{C}}\mbox{ if $\alpha=-1$}, \\ \mathfrak{N}_2^{\mathbb{C}}\left(\frac{\alpha}{(1+\alpha)^2}\right)\mbox{ otherwise}
\end{array}}&  \begin{array}{c}e_1e_3 =e_1+e_2, \  e_2e_3=\alpha e_2 \end{array}& 4\\ 
\hline

\mathfrak{g}_4&
\multicolumn{7}{|l|}{ \begin{array}{l} \mathfrak{N}_3(0)\end{array} } &
\begin{array}{c}e_1e_2 =e_3, \ e_1e_3=-e_2, \ e_2e_3=e_1 \end{array}& 3\\
\hline
\hline

\mathcal{A}_1^{\alpha}& 
\multicolumn{7}{|l|}{
\begin{array}{l}\mathfrak{N}_4\mbox{ if $\alpha=-1$},\, \mathfrak{N}_5\mbox{ if $\alpha=1$},\\ \mathfrak{N}_3\left(\frac{1+\alpha}{1-\alpha}i\right)\mbox{ otherwise}\end{array} } & \begin{array}{c}e_1e_2=e_3, \ e_1e_3 =e_1+e_3, \  e_2e_3=\alpha e_2
\end{array} & 1\\ 
\hline
\mathcal{A}_2& \multicolumn{7}{|l|}{ \begin{array}{l} \mathfrak{N}_0 \end{array} } &  \begin{array}{c}e_1e_2=e_1, \ e_2e_3=e_2\end{array} & 2\\
\hline
\mathcal{A}_3&\multicolumn{7}{|l|}{ \begin{array}{l} \mathfrak{N}_{10} \end{array} } & \begin{array}{c}e_1e_2=e_3, \ e_1e_3=e_1, \ e_2e_3=e_2\end{array} & 3\\
\hline

\multicolumn{10}{l}{}\\
\multicolumn{10}{l}{ \mbox{{\bf Table A.2.} 
{\it Classification of three dimensional  Leibniz (non-Lie) algebras.}}}\\
\multicolumn{10}{l}{}\\

\hline
\mathbb{A}& \mbox{\cite{casas}}& \multicolumn{6}{c|}{S_{\mu}} &\mbox{ Multiplication tables } & \mathfrak{Der}  \\ 
 \hline 

\mathfrak{L}^{\beta}_1 & 2(a) 
 &&&&&&
&\begin{array}{c}e_2e_2 =\beta e_1, e_3e_2=e_1, e_3e_3 =e_1\end{array} &  4
\\
\hline
\mathfrak{L}_2 & 2(b)  &&&&&&
&\begin{array}{c}e_3e_3=e_1\end{array}& 5\\
\hline
\mathfrak{L}_3 & 2(c)  &&&&&&
&\begin{array}{c}e_2e_2=e_1, e_3e_3=e_1\end{array}& 4\\
\hline

\mathfrak{L}^{\alpha}_4 & 2(e),2(f) & (&\beta,&0,&1,&-1&)
&\begin{array}{c}e_1e_3=\alpha e_1, e_2e_3=e_2, e_3e_2=-e_2,e_3e_3=e_1 \end{array}& 3 \\
\hline
\mathfrak{L}_5 & 2(g)
& (&2,&0,&1,&-1&)
&\begin{array}{c}e_1e_3=2e_1, e_2e_2=e_1, e_2e_3=e_2, e_3e_2=-e_2, e_3e_3=e_1 \end{array}& 2\\
\hline
\mathfrak{L}^{\alpha}_6 & 2(d),3(a)& (&\alpha,& 0,&1,&0&)
 &\begin{array}{c}e_1e_3=\alpha e_1, e_2e_3=e_2 \end{array}& 2+\delta_{\alpha,0}+2\delta_{\alpha,1} \\
\hline
\mathfrak{L}_7 &3(b) & (&1,&0,&1,&0&)
 & \begin{array}{c}e_1e_3=e_1+e_2, e_2e_3=e_2 \end{array}& 2\\
\hline
\mathfrak{L}_{8} &3(c) &&&&&& 
&\begin{array}{c}e_1e_3=e_2, e_3e_3=e_1 \end{array}& 3\\
\hline
\mathfrak{L}_{9} &3(d)  & (&0,&0,&1,&0&)
& \begin{array}{c}e_1e_3=e_2, e_2e_3=e_2, e_3e_3=e_1 \end{array}& 2\\
\hline

\multicolumn{10}{l}{}\\
\multicolumn{10}{l}{ \mbox{{\bf Table A.3.} 
{\it Degenerations of Leibniz algebras of dimension $3$.}}}\\

\multicolumn{10}{l}{}\\
\hline

\multicolumn{10}{|c|}{
\begin{array}{|rcl|lll|}

\multicolumn{3}{|c|}{\mbox{Degenerations}} &   \multicolumn{3}{|c|}{\mbox{Parametrized bases}}\\
\hline

\mathfrak{L}^{\beta}_1  & \to &    \mathfrak{L}_{2} &    
E_1^t=e_1, & E_2^t=te_2, & E_3^t=e_3  \\
\hline

\mathfrak{L}^{1/4}_1 & \to &  \mathfrak{g}_{1}    & 
E_1^t=t^3e_1, & E_2^t=-2te_2+te_3, & E_3^t=2t^2e_2  \\
\hline

\mathfrak{L}_3 & \to &  \mathfrak{L}_{2}    & 
E_1^t=e_1, & E_2^t=te_2, & E_3^t=e_3 \\
\hline

\mathfrak{L}^{\alpha\neq 2}_4 & \to &  \mathfrak{L}_{1}^{\frac{1-\alpha}{(2-\alpha)^2}}  

& 
E_1^t= t^2e_1, & E_2^t=\frac{t}{2-\alpha}e_1+\frac{(1-\alpha+t)t}{2-\alpha}e_3, & 
E_3^t= \frac{1}{2-\alpha}e_2+te_3 \\

\hline

\mathfrak{L}^{0}_4 & \to &  \mathfrak{g}^{0}_{3}    & 
E_1^t=-t^{-1}e_1+te_2, & E_2^t=t^{-1}e_1, & E_3^t=e_3 \\

\hline


\mathfrak{L}_4^{2} & \to &  \mathfrak{L}_{3}    & 
E_1^t=t^2e_1, & E_2^t=ite_1+e_2-ite_3, & E_3^t=te_3 \\

\hline




\mathfrak{L}_{5} & \to &  \mathfrak{L}^{2}_{4}    & 
E_1^t=te_1, & E_2^t=te_2, & E_3^t=\frac{t-1}{2}e_1+e_3 \\

\hline

\mathfrak{L}^{0}_6 & \to &  \mathfrak{L}^0_{1}    & 
E_1^t=t^2e_2, & E_2^t=te_3, & E_3^t=te_1+te_2+te_3  \\

\hline


\mathfrak{L}^{1}_6 & \to &  \mathfrak{L}_{2}    & 
E_1^t=te_1, & E_2^t=e_2, & E_3^t=e_1+te_3 \\

\hline
\mathfrak{L}^{\alpha\neq0,1}_6 & \to &  \mathfrak{L}_{8}    & 
E_1^t=e_1+te_2, & E_2^t=(\alpha-1)te_1, & E_3^t=\alpha^{-1}t^{-1}e_1+e_2+te_3 \\

\hline

\mathfrak{L}_{7} & \to &  \mathfrak{L}^{1}_{6}    & 
E_1^t=te_1, & E_2^t=e_2, & E_3^t=e_3 \\
\hline

\mathfrak{L}_{7} & \to &  \mathfrak{L}_{8}    & 
E_1^t=e_1+e_2, & E_2^t=te_2, & E_3^t=t^{-1}e_1+te_3 \\

\hline

\mathfrak{L}_{8} & \to &  \mathfrak{L}^{0}_{1}    & 
E_1^t=te_2, & E_2^t=te_3, & E_3^t=e_1+te_3 \\

\hline


\mathfrak{L}_{9} & \to &  \mathfrak{L}^0_{6}    & 
E_1^t=t^{-1}e_1, & E_2^t=t^{-2}e_2, & E_3^t=e_3 \\

\hline

\mathfrak{L}_{9} & \to &  \mathfrak{L}_{8}    & 
E_1^t=t^2e_1, & E_2^t=t^3e_2, & E_3^t=te_3 \\
\hline 
\end{array}
}\\

\hline

\end{array}$$
\end{center}

\bigskip

\bigskip

\begin{center} \footnotesize

$$\begin{array}{|lcl|llll|l|}

\multicolumn{8}{l}{}\\
\multicolumn{8}{l}{  \mbox{{\bf Table A.4.}  {\it Non-degenerations of Leibniz algebras of dimension $3$.}}} \\

\multicolumn{8}{l}{}\\

\hline
\multicolumn{3}{|c|}{\mbox{Non-degenerations}}  &  \multicolumn{5}{|c|}{\mbox{Reasons}} \\
\hline



\mathfrak{L}_4^{\alpha\neq2} & \nrightarrow & 
\mathfrak{B}, \mathfrak{B}\in \left\{
\begin{array}{l}\mathfrak{L}_1^{\beta \neq \frac{1-\alpha}{(\alpha-2)^2}},\mathfrak{L}_3, \\ \mathfrak{g}_1(\alpha\not=0), \mathfrak{g}_2,  \mathfrak{g}_3^{\beta\neq 0} \end{array}\right\} &
\multicolumn{5}{|l|}{\begin{array}{ll}
\mathcal{R}=&\left\{A\left| \begin{array}{l}
A=\langle f_1,f_2,f_3 \rangle, 
c_{11}^2=0, c_{21}^2=-c_{12}^2,\\ c_{31}^3=-\alpha c_{12}^2, c_{21}^3=(\alpha-1)c_{12}^3, \\
A_1A_3+A_2A_2=0,  A_3A_1\subseteq A_3, A_1A_1 \subseteq A_2  
\\ \end{array}\right.\right\}\vspace{0.1cm}\\
&\mathfrak{L}_4^{\alpha} \in O(\mathcal{R}) \mbox{ (take $f_1=e_3, f_2=e_2, f_3=e_1$), but } 
\mathfrak{B}  \not\in O(\mathcal{R})
\end{array}}\\


\hline

\hline

\mathfrak{L}_5 & \nrightarrow & 
\mathfrak{B}, \mathfrak{B}\in\{\mathfrak{L}_1^{\beta}, \mathfrak{g}_1, \mathfrak{g}_2\} &
\multicolumn{5}{|l|}{\begin{array}{ll}
\mathcal{R}=&\left\{A\left| \begin{array}{l}
A=\langle f_1,f_2,f_3 \rangle, 
c_{31}^3=2c_{21}^2=-2c_{12}^2, c_{21}^3=c_{12}^3, \\
A_1A_3+A_3A_2=0,  A_3A_1+A_2A_2 \subseteq A_3,\\
A_1A_2+A_2A_1 \subseteq A_2  \\ \end{array}\right.\right\}\vspace{0.1cm}\\
&\mathfrak{L}_5 \in O(\mathcal{R}) \mbox{ (take $f_1=e_3, f_2=e_2, f_3=e_1$), but } 
\mathfrak{B} \not\in O(\mathcal{R})
\end{array}}\\

\hline


\hline


\mathfrak{A} & \nrightarrow & \mathfrak{B},\begin{array}{l}\mathfrak{A}\in \{\mathfrak{L}^{\alpha}_{6},\mathfrak{L}_{7},\mathfrak{L}_{9}\},\\
\mathfrak{B}\in \left\{ 
\begin{array}{l}
\mathfrak{L}^{\beta \neq0}_1,\mathfrak{L}_3, \\ 
\mathfrak{L}_4^{1}, \mathfrak{g}_1,\mathfrak{g}_2  \end{array}\right\}\end{array} & 
\multicolumn{5}{|l|}{Ann_L(\mathfrak{A})>Ann_L(\mathfrak{B})}\\
\hline





\multicolumn{8}{l}{}\\

\multicolumn{8}{l}{ \mbox{ {\bf Table A.5.}   
{\it Orbit closures for some families of three dimensional Leibniz algebras.}}} \\

\multicolumn{8}{l}{}\\

\hline
\multicolumn{3}{|c|}{\mbox{Degenerations}}  &  \multicolumn{4}{c|}{\mbox{Parametrized bases}} &  \multicolumn{1}{c|}{\mbox{Indices}} \\

\hline

\mathfrak{L}_4^{*} &\to&  \mathfrak{L}^0_{6} & 
E_1^t=e_2, & E_2^t=e_1, & E_3^t=te_3 && \epsilon=t^{-1}\\


\hline
\mathfrak{L}_6^{*} &\to&  \mathfrak{L}_{7} & 
E_1^t=e_1+e_2, & E_2^t=te_2, & E_3^t=e_3 && \epsilon=1-t\\
\hline

\mathfrak{L}_6^{*} &\to&  \mathfrak{L}_{9} & 
E_1^t=e_1+te_2, & E_2^t=(1-t)e_1, & E_3^t=e_1+e_2+te_3 && \epsilon=t^{-1}\\
\hline

\mathfrak{L}_1^{*} &\to&  \mathfrak{L}_{3} & 
E_1^t=t^4e_1, & E_2^t=t^3e_2, & E_3^t=t^2e_3 && \epsilon= t^{-2}\\
\hline

\hline
\hline
\multicolumn{3}{|c|}{\mbox{Non-degenerations}}  &  \multicolumn{5}{c|}{\mbox{Reasons}} \\
\hline

\hline





\mathfrak{L}^{*}_{4} & \nrightarrow &\mathfrak{B}, \mathfrak{B}\in\{\mathfrak{L}_{6}^1, \mathfrak{L}_{8}\}   & 
\multicolumn{5}{l|}{(\mathfrak{L}^{\alpha}_{4})^{(+2)}<\mathfrak{B} ^{(+2)}}\\
\hline
\mathfrak{L}^{*}_{4} & \nrightarrow &  \mathfrak{B}, \mathfrak{B}\in\{\mathfrak{g}_2,\mathfrak{g}^{\beta\not=0}_3\}   & 
\multicolumn{5}{l|}{Ann_L(\mathfrak{L}^{\alpha}_{4})>Ann_L(\mathfrak{B})}\\
\hline

\mathfrak{L}^{*}_{6} & \nrightarrow &  \mathfrak{B}, \mathfrak{B}\in\{\mathfrak{g}_1,\mathfrak{g}_2,\mathfrak{L}_{1}^{\beta\not=0},\mathfrak{L}_{3},\mathfrak{L}_{4}^{1}\}   & 
\multicolumn{5}{l|}{Ann_L(\mathfrak{L}^{\alpha}_{6})>Ann_L(\mathfrak{B})}\\
\hline

\end{array}$$
\end{center}

\end{document}